\documentclass[12pt]{article}

\usepackage{geometry} 
\usepackage{latexsym}
\usepackage{amssymb}
\usepackage{amsthm}
\usepackage{amsmath,amsfonts}
\usepackage{xcolor}
\usepackage{float}

\def\titlerunning#1{\gdef\titrun{#1}}
\makeatletter
\def\author#1{\gdef\autrun{\def\and{\unskip, }#1}\gdef\@author{#1}}
\def\address#1{{\def\and{\\\hspace*{18pt}}\renewcommand{\thefootnote}{}%
\footnote {#1}}%
\markboth{\autrun}{\titrun}}
\makeatother
\def\email#1{\hspace*{4pt}{\em e-mail}: #1}

\titlerunning{}

\title{New directed strongly regular graphs on 60 vertices}
\author{Dean Crnkovi\' c, Andrea \v Svob and Matea Zubovi\'c \v Zutolija}

\begin{document}
\maketitle

\address{D. Crnkovi\' c, A. \v Svob, M. Zubovi\'c \v Zutolija: Faculty of Mathematics, University of Rijeka, Radmile Matej\v ci\'c 2, 51000 Rijeka, Croatia;
\email{\{deanc,asvob,matea.zubovic\}@math.uniri.hr}
}

\begin{abstract}
We prove the existence of directed strongly regular graphs with parameters (60,21,11,6,8), (60,22,12,8,8), (60,24,10,9,10), (60,25,17,8,12), (60,27,21,12,12) and (60,28,20,14,12). The group $S_5 \times 2$ acts transitively on the constructed graphs.
\end{abstract}

\bigskip



\section{Introduction}

A directed strongly regular graph (dsrg) with parameters $v,k,t,\lambda,\mu$ is a directed graph on $v$ vertices, such that each vertex has indegree and outdegree $k$,
and for any two vertices $x,y$ the number of directed paths $x \to z \to y$ is $t$, $\lambda$ or $\mu$ when $x=y$, or $xy$ is an arc, or $xy$ is not an arc, respectively. The digraphs that we use in this paper will not have more than one arc from one vertex to another, and further, will not have any arcs from a vertex to itself.

Directed strongly regular graphs were introduced by Duval \cite{Duval88}, as a directed version of strongly regular graphs.
According to the literature and the table in \cite{BH}, no directed strongly regular graphs with parameters (60,21,11,6,8), (60,22,12,8,8), (60,24,10,9,10), (60,25,17,8,12), (60,27,21,12,12) and (60,28,20,14,12) were known. In this note, we construct five nonisomorphic directed strongly regular graphs with parameters (60,21,11,6,8), 26 nonisomorphic directed strongly regular graphs with parameters (60,22,12,8,8), 32 nonisomorphic directed strongly regular graphs with parameters (60,24,10,9,10), four nonisomorphic directed strongly regular graphs with parameters (60,25,17,8,12), 24 nonisomorphic directed strongly regular graphs with parameters (60,27,21,12,12) and six nonisomorphic directed strongly regular graphs with parameters (60,28,20,14,12) that are the first known directed strongly regular graphs with these parameters.

\bigskip

The computations in this paper are made by using programs written for Magma \cite{magma}.

The directed strongly regular graphs constructed in this paper can be found at the link:

\begin{verbatim}
https://www.math.uniri.hr/~asvob/DSRGs_v60.html
\end{verbatim}

\section{Construction of new directed strongly regular graphs}

We constructed the directed strongly regular graphs by using the method described in \cite[Theorem 3]{cms}.
That construction produces simple 1-designs on which a group $G$ acts transitively on the points and blocks. Hence, if the incidence structure of a 1-design obtained using \cite[Theorem 3]{cms} is the adjacency matrix of a directed strongly regular graph, then the graph admits a transitive action of $G$ on the set of vertices. In this paper, we construct the directed strongly regular graphs using the group $G=S_5 \times 2$.

\subsection{Construction of dsrg(60,28,20,14,12)}

The group $S_5 \times 2$ has seven conjugacy classes of subgroups isomorphic to $2^2$. Up to conjugation, there is exactly one subgroup $H_1$ isomorphic to $2^2$ having the property that $S_5 \times 2$ acts on cosets of $H_1$ in 19 orbits. Using the method from \cite[Theorem 3]{cms}, by taking $G=S_5 \times 2$ and the stabilizer of a vertex $G_{\alpha}=H_1$, we constructed two nonisomorphic 
dsrg(60,28,20,14,12) admitting a transitive action of the group $S_5 \times 2$, which we denote by $\Delta_1$ and $\Delta_2$, where one digraph is obtained from the other by reversing the arcs. We say that the digraphs $\Delta_1$ and $\Delta_2$ are reverse to each other (see \cite[Section 3.4]{BCSZ}). The full automorphism group of the digraphs $\Delta_1$ and $\Delta_2$ is isomorphic to $S_5 \times 2$, and $A_5$ acts regularly on the vertices of the digraphs.

The group $G$ has, up to conjugation, exactly one subgroup $H_2$ isomorphic to $2^2$ having the property that $S_5 \times 2$ acts on cosets of $H_2$ in 32 orbits. Using the method from \cite[Theorem 3]{cms}, by taking $G=S_5 \times 2$ and the stabilizer of a vertex $G_{\alpha}=H_2$, we constructed four nonisomorphic 
dsrg(60,28,20,14,12) admitting a transitive action of the group $S_5 \times 2$, which we denote by $\Delta_3$, $\Delta_4$, $\Delta_5$ and $\Delta_6$. The digraph $\Delta_4$ can be obtained by reversing the arcs of $\Delta_3$, and vice versa. The same holds for digraphs $\Delta_5$ and $\Delta_6$. 
The full automorphism group of the digraphs $\Delta_3$ and $\Delta_4$ is isomorphic to $S_5 \times 2$, and $A_5$ acts in two orbits on the vertices of the digraphs. The full automorphism group of the digraphs $\Delta_5$ and $\Delta_6$ is isomorphic to $S_5$, and $A_5$ acts in two orbits on the vertices of the digraphs.

The group $G$ has, up to conjugation, two subgroups isomorphic to $2^2$ such that $G$ acts on its cosets in 22 orbits. Let us denote these subgroups with $H_3$ and $H_4$.  
We obtained another pair of directed strongly regular graphs with parameters (60,28,20,14,12), isomorphic to the ones denoted by $\Delta_3$ and $\Delta_4$, by taking the subgroup $H_3$ as the stabilizer of a vertex. The subgroup $H_4$ produces two directed strongly regular graphs with parameters (60,28,20,14,12) that are isomorphic to $\Delta_3$ and $\Delta_4$.

\subsection{Construction of dsrg(60,22,12,8,8)}

By taking the subgroup $H_2$ as the stabilizer of a vertex, we obtained 26 nonisomorphic 
dsrg(60,22,12,8,8) admitting a transitive action of the group $S_5 \times 2$, which we denote by $\Delta_7, \dots, \Delta_{32}$.

The full automorphism group of the digraphs $\Delta_{17}$, $\Delta_{18}$, $\Delta_{19}$ and $\Delta_{20}$  is isomorphic to $S_5 \times 2$, while the full automorphism group for the rest of the digraphs is isomorphic to $S_5$. The group $A_5$ acts in two orbits on the vertices of the digraphs. Among these 26 nonisomorphic digraphs, there are 13 pairs of reversed digraphs.

By taking any of the subgroups $H_3$ or $H_4$ as the stabilizer of a vertex, we obtained two pairs of pairwise reversed directed strongly regular graphs with parameters (60,22,12,8,8), isomorphic to digraphs $\Delta_{17}$, $\Delta_{18}$, $\Delta_{19}$ and $\Delta_{20}$.

\subsection{Construction of dsrg(60,25,17,8,12)}

By applying \cite[Theorem 3]{cms} to $G=S_5 \times 2$, taking the subgroup $H_2$ as the stabilizer of a vertex, we also obtained two pairs of reversed directed strongly regular graphs with parameters (60,25,17,8,12), which we denote with $\Delta_{33}$ and $\Delta_{34}$, and $\Delta_{35}$ and $\Delta_{36}$. The full automorphism group of $\Delta_{33}$ and $\Delta_{34}$ is $S_5$, and the full automorphism group of digraphs $\Delta_{35}$ and $\Delta_{36}$ is $G$, and $A_5$ acts on the vertices of these digraphs in two orbits.

By applying \cite[Theorem 3]{cms} to $G=S_5 \times 2$, taking any of the subgroups $H_3$ or $H_4$ as the stabilizer of a vertex, we also obtained a pair of reversed directed strongly regular graphs with parameters (60,25,17,8,12), isomorphic to digraphs $\Delta_{33}$ and $\Delta_{34}$.

\subsection{Construction of dsrg(60,24,10,9,10)}

By taking the subgroup $H_2$ as the stabilizer of a vertex, we obtained 32 nonisomorphic 
dsrg(60,24,10,9,10) admitting a transitive action of the group $S_5 \times 2$, which we denote by $\Delta_{37}, \dots, \Delta_{68}$. 

The full automorphism group of the digraphs is isomorphic to $S_5$, and the group $A_5$ acts in two orbits on the vertices of the digraphs. Among these 32 nonisomorphic digraphs, there are 16 pairs of reversed digraphs.

\subsection{Construction of dsrg(60,27,21,12,12)}

By taking the subgroup $H_2$ as the stabilizer of a vertex, we obtained 24 nonisomorphic 
dsrg(60,27,21,12,12) admitting a transitive action of the group $S_5 \times 2$, which we denote by $\Delta_{69}, \dots, \Delta_{92}$.

The full automorphism group of the digraphs is isomorphic to $S_5$, and the group $A_5$ acts in two orbits on the vertices of the digraphs. Among these 24 nonisomorphic digraphs, there are 12 pairs of reversed digraphs.

\subsection{Construction of dsrg(60,21,11,6,8)}

By taking the subgroup $H_2$ as the stabilizer of a vertex, we obtained five nonisomorphic 
dsrg(60,21,11,6,8) admitting a transitive action of the group $S_5 \times 2$, which we denote by $\Delta_{93}, \dots, \Delta_{97}$.

The full automorphism group of the digraphs is isomorphic to $S_5$, and the group $A_5$ acts in two orbits on the vertices of the digraphs. Among these five nonisomorphic digraphs, we obtained two pairs of reversed digraphs, while one digraph is isomorphic to its reverse.

\section{Summarizing results}

We can conclude that in all cases where there is an even number of nonisomorphic graphs, none of them is isomorphic to its reverse, while in all cases where there is an odd number of nonisomorphic graphs,
precisely one example is isomorphic to its reverse. Note that the only open cases for directed strongly regular graphs on 60 vertices are now $(60,16,9,2,5)$, $(60,17,11,4,5)$, $(60,18,11,6,5)$ and $(60,22,19,6,9)$, since with the results from this paper, we covered six open cases for parameter sets. In Table \ref{new}, we summarize the results. 

\begin{table}[H]
\centering
\begin{tabular}{|c|c|c|c|c|}
\hline
$H$ & $(n,k,t,\lambda,\mu)$ & $\#$ nonisom. & Aut$\,\Gamma$\\
	\hline
$H_2$ & (60,21,11,6,8) & 5 & $S_5$  \\
$H_2$ & (60,22,12,8,8) & 26 & $S_5\times 2$ (4); $S_5$ (22)  \\
$H_2$ & (60,24,10,9,10) & 32 & $S_5$ \\
$H_2$ & (60,25,17,8,12) & 4 & $S_5\times 2$ (2); $S_5$ (2) \\
$H_2$ & (60,27,21,12,12) & 24 & $S_5$  \\
$H_2$ & (60,28,20,14,12) & 4 & $S_5\times 2$ (2); $S_5$ (2)  \\
$H_1$ & (60,28,20,14,12) & 2 & $S_5\times 2$   \\
\hline
\end{tabular}
\caption{New DSRGs parameter sets from $S_5 \times 2$}\label{new}
\end{table}

\section*{Acknowledgement}
This work was supported by the Croatian Science Foundation under the project number HRZZ-IP-2022-10-4571 and by European Union-NextGenerationEU, project number uniri-iz-25-46-KonGeoGraGru.

\end{document}